%
%
%




\input amstex.tex

\magnification=\magstep1
\hsize=5.5truein
\vsize=9truein
\hoffset=0.5truein
\parindent=10pt
\newdimen\nagykoz
\newdimen\kiskoz
\nagykoz=7pt
\kiskoz=2pt
\parskip=\nagykoz
\baselineskip=12.7pt


\loadeufm \loadmsam \loadmsbm

\font\vastag=cmssbx10
\font\drot=cmssdc10
\font\vekony=cmss10
\font\vekonydolt=cmssi10
\font\cimbetu=cmssbx10 scaled \magstep1
\font\szerzobetu=cmss10

\font\scVIII=cmcsc8
\font\rmVIII=cmr8
\font\itVIII=cmti8
\font\bfVIII=cmbx8
\font\ttVIII=cmtt8

\def\cim#1{{\centerline{\cimbetu#1}}}
\def\szerzo#1{{\vskip0.3truein\centerline{\szerzobetu#1}}}
\def\alcim#1{{\medskip\centerline{\vastag#1}}}
\def\tetel#1#2{{{\drot#1}{\it\szukebb~#2\tagabb}}}
\long\def\biz#1#2{{{\vekony#1} #2}}
\def\kiemel#1{{\vekonydolt#1\/}}
\long\def\absztrakt#1#2{{\vskip0.4truein{\vekony#1} #2\vskip0.5truein}}
\def\szukebb{\parskip=\kiskoz}
\def\tagabb{\parskip=\nagykoz}
\def\vonal{{\vrule height 0.2pt depth 0.2pt width 0.5truein}}

\def\CC{{\Bbb C}}

\def\ts{\textstyle}

\def\bmfd{{Banach manifold}}

\def\bvbdl{{Banach vector bundle}}

\def\Oa{{\Omega}}

\def\UU{{\frak U}}

\def\BB{{\frak B}}

\def\Hom{\hbox{\rm Hom}}
\def\End{\hbox{\rm End}}

\def\cts{{continuous}}
\def\idml{{infinite dimensional}}

\def\fdml{{finite dimensional}}

\def\pscx{{pseudoconvex}}

\def\dist{{\hbox{\rm dist}}}

\def\sbs{{Schauder basis}}

\def\cubs{{countable unconditional basis}}
\def\psh{{plurisubharmonic}}
\def\pshdom{{\psh\ domination}}

\def\st{{such that}}
\def\wrt{{with respect to}}

\def\<{{\langle}}
\def\>{{\rangle}}

\def\RR {{\Bbb R}}

\def\AA{{\Cal A}}

\def\OO {{\Cal O}}
\def\LL {{\Cal L}}

\def\ro{\varrho}

\def\epsz{\varepsilon}
\def\fii{\varphi}
\def\fn{func\-tion}
\def\fns{func\-tions}
\def\holo{hol\-o\-mor\-phic}

\def\nbd{neighbor\-hood}

\def\bspc{Banach space}

\def\la{\lambda}
\def\La{\Lambda}

\def\da{{\delta}}

\def\ga{{\gamma}}
\def\aa{{\alpha}}
\def\ba{{\beta}}
\def\da{{\delta}}

\def\Prop{Proposition}
\def\p#1{{\Prop~#1}}

\def\Th{{Theorem}}

\def\t#1{{\Th~#1}}

{\phantom.}
\vskip0.5truein
\cim{AN ANALYTIC KOSZUL COMPLEX}
\vskip0.2truein
\cim{IN A BANACH SPACE}
\szerzo{Imre Patyi\plainfootnote{${}^*$}{\rmVIII Supported in part by a Research Initiation Grant
	from Georgia State University.}}
\absztrakt{ABSTRACT.}{We show that the holomorphic ideal sheaf of a linear
	section of a pseudoconvex open subset $\Omega$ of, say, a Hilbert space
	$X=\ell_2$ is acyclic.
	 We also prove an analog of Hefer's lemma, i.e., if $f:\Omega\times\Omega\to\CC$
	is holomorphic and $f(x,x)=0$ for $x\in\Omega$, then there is a holomorphic 
	$g:\Omega\times\Omega\to X^*$ with values in the dual space $X^*$ of $X$ such that
	$f(x,y)=g(x,y)(x-y)$.

	
}

\def\sA{{1}}
\def\sB{{2}}
\def\sC{{3}}
\def\sD{{4}}
\def\sE{{5}}


\def\tAA{{\sA.1}}
\def\tAB{{\sA.2}}

\def\eAA{{\sA.1}}

\def\tBA{{\sB.1}}
\def\tBB{{\sB.2}}
\def\tBC{{\sB.3}}
\def\tBD{{\sB.4}}
\def\tBE{{\sB.5}}

\def\eBA{{\sB.1}}
\def\eBB{{\sB.2}}
\def\eBC{{\sB.3}}

\def\tCA{{\sC.1}}
\def\tCB{{\sC.2}}

\def\eCA{{\sC.1}}
\def\eCB{{\sC.2}}
\def\eCC{{\sC.3}}
\def\eCD{{\sC.4}}

\def\tDA{{\sD.1}}
\def\tDB{{\sD.2}}
\def\tDC{{\sD.3}}

\def\eDA{{\sD.1}}
\def\eDB{{\sD.2}}
\def\eDC{{\sD.3}}

\def\tEA{{\sE.1}}
\def\tEB{{\sE.2}}
\def\tEC{{\sE.3}}
\def\tED{{\sE.4}}
\def\tEE{{\sE.5}}

\def\eEA{{\sE.1}}
\def\eEB{{\sE.2}}

\def\rLA{L1}
\def\rLB{L2}
\def\rLC{L3}
\def\rNW{NW}
\def\rPA{P1}
\def\rPB{P2}

\alcim{\sA. INTRODUCTION.}

	 Variants of the Koszul complex play so important a role
	in commutative algebra, and algebraic and analytic geometry that
	one can verily call it the mother of all resolutions.

	 This paper establishes an exactness and acyclicity result for
	an analytic Koszul complex in a \bspc\ that serves as preparation
	to show in [\rPB] that the ideal sheaf of certain complex submanifolds
	of a \bspc\ belongs to a class of sheaves, to be called therein of type~(S),
	which are studied and proved acyclic in certain cases therein via a
	method of resolutions.
	 In effect, we do here the case of the ideal sheaf of a linear submanifold.

	 Let $X',X'',Z$ be complex \bspc{}s, $X=X'\times X''$,
	$\La_p$ the \bspc\ of all \cts\ complex $p$-linear
	alternating maps $X''\to Z$ for $p\ge0$; $\La_0=\La_{-1}=Z$; and $\OO^{\La_p}\to X$ the
	sheaf of germs of \holo\ functions $X\to\La_p$.
	 Let $E$ be the Euler vector 
	field on $X''$ defined by $E(x'')=x''$, 
	$\LL_E$ the Lie derivation,  and
	$i_E$ the inner
	derivation determined by the vector field $E$, i.e., $i_E$ is the contraction
	of $p$-forms with $E$: if $f$ is a local section of $\OO^{\La_p}$, then
	let $i_Ef$ be the local section of $\OO^{\La_{p-1}}$ given for $p\ge1$ by
	$(i_Ef)(x',x'')(\xi''_1,\xi''_2,\ldots,\xi''_{p-1})=
	 f(x',x'')(x'',\xi''_1,\ldots,\xi''_{p-1})$, and for $p=0$ by
	$(i_Ef)(x',x'')=f(x',0)$.
	 Let $I$ be the subsheaf of $\OO^Z$ of all sections that vanish on $X'$.
	 We consider the Koszul complex
$$
	\ldots\to\OO^{\La_p}\to\OO^{\La_{p-1}}\to\ldots\to\OO^{\La_1}\to I\to0
\tag\eAA
 $$
 	of analytic sheaves over $X$, where each map is $i_E$.
	 Let $K_p$, $p\ge0$, be the corresponding sequence of kernel sheaves:
	$K_p(U)=\{f\in\OO(U,\La_p):i_Ef=0{\text{\/ on\ }}U\}$, $U\subset X$ open; $K_0=I$.

	 Lempert [\rLA] introduced the notion of \pshdom, and demonstrated its
	usefulness for proving vanishing theorems first in [\rLB].
	 Following him we say that \kiemel{\pshdom} holds in a complex \bmfd\ $\Oa$
	if given any locally upper bounded \fn\ $u:\Oa\to\RR$ there is a \cts\
	\psh\ \fn\ $\psi:\Oa\to\RR$ \st\ $u(x)<\psi(x)$ for all $x\in\Oa$.

\tetel{\t\tAA.}{{\rm (Lempert)}
\vskip0pt
{\rm (a) [\rLA]} If $X$ is a \bspc\ with a \cubs, and $\Oa\subset X$ is \pscx\ open,
	then \pshdom\ holds in\/ $\Oa$.
\vskip0pt
{\rm (b) [\rLB]} If $X$ is a \bspc\ with a \sbs, \pshdom\ holds in an open subset\/ $\Oa$ of $X$
	(consequently, $\Oa$ is \pscx), and $Z$ is any \bspc, then
	the sheaf cohomology groups $H^q(\Oa,\OO^Z)$ vanish for all $q\ge1$, where
	$\OO^Z\to\Oa$ is the sheaf of germs of \holo\ \fn{}s\/ $\Oa\to Z$.
}

	(Lempert proved in [\rLB] a theorem different from \t\tAA(b), but his
	method can be modified to give \t\tAA(b) as it stands above; see [\rPA, Thm.\,1.3(a)].)
	 
	 \t\tAB\ below is this paper's main theorem.

\tetel{\t\tAB.}{With the above notation, let\/ $\Oa\subset X$ be \pscx\
	open.  If $X',X''$ have Schauder bases, and \pshdom\ holds in every \pscx\ open
	subset of\/ $\Oa$, then
{\vskip0pt\rm (a)}
	the Koszul complex\/ $(\eAA)$ is exact on the germ level and on the level
	of global sections over\/ $\Oa$, and
{\vskip0pt\rm (b)}
	the $K_p$ are acyclic over\/ $\Oa$: $H^q(\Oa,K_p)=0$ for all $q\ge1$ and $p\ge0$.
}

	 Note that \t\tAA\ implies that \t\tAB\ applies
	when $X',X''$ have countable unconditional bases; e.g., when $X=\ell_2$.
	 In the proof of \t\tAB\ we first look at the local exactness of $(\eAA)$,
	then prove vanishing by a Leray covering argument combined with Lempert's
	method of exhaustion.
	 See \t\tEA\ for an extension of \t\tAB.

\alcim{\sB. LOCAL EXACTNESS AND VANISHING.}

	 In this section we show that the Koszul complex $(\eAA)$ is exact on the
	germ level and on the level of global sections over suitable \pscx\ open
	\nbd{}s $\Oa'$ of each point $x_0\in\Oa$.

\tetel{\p\tBA.}{Let\/ $\Oa'\subset\Oa$ be \pscx\ open.
	 If\/ {\rm(a)} or\/ {\rm(b)} below holds, then\/ $(\eAA)$ is exact
	 on the level of global sections over\/ $\Oa'$.
\vskip0pt
	{\rm(a)} $\Oa'$ is such that there is a $\la\in\OO(X'')$
	that is linear, or more generally homogeneous of degree $m\ge1$,
	with $\la(x'')\not=0$ for $(x',x'')\in\Oa'$.
\vskip0pt
	{\rm(b)} $\Oa'$ is such that the flow lines of the Euler vector field $E$
	stay in\/ $\Oa'$ for time $-\infty\le t\le0$, i.e., if\/ $(x',x'')\in\Oa'$,
	then\/ $(x',tx'')\in\Oa'$ for all\/ $0\le t\le1$.
}

\biz{Proof.}{Given any $f\in\OO(\Oa',\La_p)$ with $i_Ef=0$ on $\Oa'$ we need to produce
	a $g\in\OO(\Oa',\La_{p+1})$ with $i_Eg=f$ on $\Oa'$.  
	Let $d=d_{x''}$ be the usual operator of outer differentiation.

	(a) Letting $g=\frac{d\la}{m\la}\wedge f$ will do.
	 Here 
$$\eqalign{
	g(x',x'')(\xi_0'',\xi_1'',&\ldots,\xi_p'')=\cr
	& \sum_{i=0}^p(-1)^i\frac{d\la}{m\la}(x'')(\xi_i'')f(x',x'')(\xi_0'',\ldots,\xi_{i-1}'',
	 \xi_{i+1}'',\ldots,\xi_p'').\cr
}
 $$
 	Then $i_Eg=f$ since
	$(i_Eg)(x',x'')(\xi_1'',\ldots,\xi_p'')=\frac{\la(x'')}{\la(x'')}f(x',x'')(\xi_1'',\ldots,\xi_p'')
	+0+\ldots+0=f(x',x'')(\xi_1'',\ldots,\xi_p'')$.

	(b) If $p=0$, then $i_Edf=\LL_Ef$ since both sides are just the derivative of the \fn\
	$f$ \wrt\ $E$.
	 For $p\ge1$ recall the Cartan identity
	$d(i_Ef)+i_E(df)=\LL_Ef$, which is true also in any \bspc.
	 As $i_Ef=0$, we have that $i_Edf=\LL_Ef$.
	 Our $g$ will be a suitable integral of $df$ that inverts the Lie derivation $\LL_E$.
	 Let $F^t_E$ be the flow of the Euler vector field $E$, i.e.,
	$F^t_E(x'')=e^tx''$.
	 Then $(F^t_E)^*(\LL_Ef)=\frac{d}{dt}((F^t_E)^*f)$, and
	$(F^t_E)_*E=E$, where as usual $(F^t_E)^*$ is the pull back of forms, and
	$(F^t_E)_*$ is the push forward of tangent vectors by the diffeomorphism
	(biholomorphism) $F^t_E$.
	 Recall another general identity: $(F^t_E)^*(i_Edf)=i_{(F^t_E)_*E}((F^t_E)^*df)$.

	 Define $g\in\OO(\Oa',\La_{p+1})$ by
	$g=\int_{t=-\infty}^0(F^t_E)^*(df)\,dt$.
	 This integral converges and is \holo\ since on substituting $e^tx''$ in the
	form $df$ we gain at least one factor of $e^t$ in the integrand, and letting
	$s=e^t$ we can rewrite $g$ as a proper integral over $[0,1]$ \wrt\ $s$.
	Now 
$$\eqalign{
	i_Eg&=\int_{t=-\infty}^0 i_E((F^t_E)^*df)\,dt=
	 \int_{-\infty}^0i_{(F^t_E)_*E}((F^t_E)^*df)\,dt\cr&=
	 \int_{-\infty}^0(F^t_E)^*(i_Edf)\,dt=
	 \int_{-\infty}^0(F^t_E)^*(\LL_Ef)\,dt\cr&=
	 \int_{-\infty}^0\frac{d}{dt}((F^t_E)^*f)\,dt=
	 (F^0_E)^*f-(F^{-\infty}_E)^*f=f,
}
 $$
	where $F^0_E$ and $(F^0_E)^*$ are the identity, and
	$(F^{-\infty}_E)^*f=0$ because we have that
	$((F^{-\infty}_E)^*f)(x',x'')=f(x',0)=0$ for $p=0$ since $f\in I(\Oa')$, and
	$((F^{-\infty}_E)^*f)(x',x'')(\xi_1'',\ldots,\xi_p'')=f(x',0)(0,\ldots,0)=0$
	for $p\ge1$.
	 The proof of \p\tBA\ is complete.
}

	 Note that a Koszul complex can also be built over a suitable
	\pscx\ open neighborhood of a split complex Banach submanifold $M$ of $\Oa$
	if a suitable global \holo\ vector field $E$ can be found at least on
	a \nbd\ of $M$ in $\Oa$ that can play the role of the Euler vector field.
	 This is the case, e.g., when $M$ is a complete intersection in $\Oa$.

\tetel{\p\tBB.}{The sequence\/ $(\eAA)$ is exact on the germ level over\/ $\Oa$.
}

\biz{Proof.}{It is enough to show that each point $x_0=(x_0',x_0'')\in\Oa$
	has arbitrary small \nbd{}s $\Oa'$ over which $(\eAA)$ is exact
	on the level of global sections over $\Oa'$.
	 We put a norm $\|x\|=\max\{\|x'\|,\|x''\|\}$ on $x=(x',x'')\in X=X'\times X''$
	so that the open balls $B_X(x,\epsz)$ are direct products
	$B_X(x,\epsz)=B_{X'}(x',\epsz)\times B_{X''}(x'',\epsz)$.
	 We can now choose such \nbd{}s $\Oa'$ in the form of balls
	$\Oa'=B_X(x_0,\epsz)$ as follows.

	 If $x_0''\not=0$, then there is a linear functional $\la\in(X'')^*$ \st\
	$\|\la\|=1$ and $\la(x_0'')=\|x_0''\|$.
	 Then $\la(x'')\not=0$ for $\|x''-x_0''\|<\epsz<\|x_0''\|$, since
	$|\la(x'')-\|x_0''\||=|\la(x'')-\la(x_0'')\|\le\epsz<\|x_0''\|$.
	 That is, $\Oa'=B_X(x_0,\epsz)$ for $0<\epsz<\|x_0''\|$ satisfies
	\p\tBA(a) and lies in any \nbd\ of $x_0$ in $\Oa$ if $\epsz>0$ is
	small enough.

	 If $x_0''=0$, then $\Oa'=B_X(x_0,\epsz)$, being the product of a set
	in $X'$ by a convex set in $X''$ that contains $0\in X''$,
	satisfies \p\tBA(b) and lies in
	any \nbd\ of $x_0$ in $\Oa$ if $\epsz>0$ is small enough.
	 The proof of \p\tBB\ is complete.
}

	 Consider now an exact sequence
$$
	\ldots\to\OO^{L_p}\to\OO^{L_{p-1}}\to\ldots\OO^{L_1}\to K'_0\to0
\tag\eBA
 $$
 	of analytic sheaves over $\Oa$, where each differential is called $d$, and
	$L_p$ are arbitrary \bspc{}s, and the short exact sequences
$$
	0\to K'_p\to\OO^{L_p}\to K'_{p-1}\to0
\tag\eBB
 $$
 	of analytic sheaves, where $K'_p$ is the kernel of $d:\OO^{L_p}\to\OO^{L_{p-1}}$ for $p\ge2$
	and of $d:\OO^{L_1}\to K'_0$ for $p=1$.
	 By \p\tBB\ our Koszul complex $(\eAA)$ is as $(\eBA)$.

\tetel{\p\tBC.}{Let\/ $\Oa'\subset\Oa$ be \pscx\ open.
\vskip0pt
	{\rm(a)} If\/ $(\eBA)$ is exact on the level of global sections over\/
	$\Oa'$, then we have $H^q(\Oa',K'_p)=0$ for all $q\ge1$ and $p\ge q$.
\vskip0pt
	{\rm(b)} If $H^1(\Oa',K'_p)=0$ for all $p\ge1$, then\/ $(\eBA)$ is
	exact on the level of global sections over\/ $\Oa'$.
\vskip0pt
	{\rm(c)} If $H^q(\Oa',K'_p)=0$ for all $q\ge2$ and $p\ge1$, then
	$H^q(\Oa',K'_p)=0$ for all $q\ge1$ and $p\ge0$.
}

\biz{Proof.}{Writing down the long exact sequence associated to $(\eBB)$ for $p\ge1$
	in cohomology we get that
$$\eqalign{
	0&\to K'_p(\Oa')\to\OO(\Oa',L_p)\to K'_{p-1}(\Oa')\to\cr
	&\to H^1(\Oa',K'_p)\to0\to H^1(\Oa',K'_{p-1})\to\cr
	&\to H^2(\Oa',K'_p)\to0\to H^2(\Oa',K'_{p-1})\to\cr
	&\phantom{\to}\ldots\quad\ldots\quad\ldots\quad\ldots\quad\ldots\quad\ldots\quad\ldots\cr
}\quad,
\tag\eBC
 $$
 	where the zeros except the first are by Lempert's \t\tAA(b).

	(a) As the third map in $(\eBC)$ is an epimorphism by assumption,
	we see that $H^1(\Oa',K'_p)=0$ for $p\ge1$.
	 Since $H^q(\Oa',K'_{p-1})\cong H^{q+1}(\Oa',K'_p)$ for $q\ge1$ and
	$p\ge1$, we find that $H^q(\Oa',K'_p)=0$ for $q\ge1$ and $p\ge q$
	as claimed.

	(b) As $H^1(\Oa',K_p)=0$ for $p\ge1$ by assumption, we see from $(\eBC)$
	that the third map is an epimorphism, thus $(\eBA)$ is exact on
	the level of global sections over $\Oa'$ as claimed.

	(c) This follows from the dimension shifting relation 
	$H^q(\Oa',K'_{p-1})\cong H^{q+1}(\Oa',K'_p)$ for $q\ge1$ and $p\ge1$ in (a).
	 The proof of \p\tBC\ is complete.
}

\tetel{\p\tBD.}{Let\/ $\Oa'\subset\Oa$ be as in the proof of \p\tBB.
	 Then the kernel sheaves $K_p$ of the Koszul complex\/ $(\eAA)$
	are acyclic over\/ $\Oa'$ for $p\ge0$.
}

\biz{Proof.}{By the proof of \p\tBC\ it is enough to prove that $I=K_0$ is acyclic over $\Oa'$.

	 If $\Oa'\cap X'=\emptyset$, then $I=\OO^Z$ over $\Oa'$, and we are done by \t\tAA(b).

	 If $\Oa'\cap X'\not=\emptyset$, then since $\Oa'$ is a product set
	 $\Oa'=B_{X'}(x'_0,\epsz)\times B_{X''}(x''_0,\epsz)$, any $h\in\OO(\Oa'\cap X',Z)$
	 has an automatic \holo\ extension $\tilde h\in\OO(\Oa',Z)$ defined
	 simply by $\tilde h(x',x'')=h(x')$. 
	  Hence $H^1(\Oa',I)=0$.
	  Considering the short exact sequence
	$0\to I\to{}_{\Oa'}\OO^Z\to{}_{\Oa'\cap X'}\OO^Z\to0$
	of analytic sheaves over $\Oa'$ we see that $H^q(\Oa',I)=0$ for $q\ge2$
	since $H^{\ge1}(\Oa'\cap X',\OO^Z)=0$ by \t\tAA(b).
	 The proof of \p\tBD\ is complete.
}

	 Note that if $\Oa'$ has a finite covering by \pscx\ open subsets
	that is a Leray covering for each $K'_p$ for $p\ge0$, then
	in \p\tBC(a) the range `$p\ge q$' can be replaced by `$p\ge0$'
	as then the cohomology groups $H^q(\Oa',K'_p)$ vanish for
	all high enough $q$ and all $p\ge0$.

	 Remark that a finite intersection $\Oa'=\bigcap\Oa'_i$ of balls $\Oa'_i$
	as in the proof of \p\tBB\ is again a set $\Oa'$ to which \p\tBA\
	applies.
	 That is, we have a covering of $\Oa$ by balls that is a Leray covering
	of $\Oa$ for the sheaves $K_p$, $p\ge0$, as in \p\tBD.

	 We end this section with a vanishing result in the midrange between local
	and global.

\tetel{\p\tBE.}{Let\/ $\Oa'\subset\Oa$ be \pscx\ open.
	Suppose that there is a finite set\/ $\UU$ of \pscx\ open subsets $U$ of\/
	$\Oa$ that is a Leray covering of $\bigcup\UU\supset\Oa'$ for $K'_p$, $p\ge1$,
	as in\/ $(\eBB)$.
	 Then $H^q(\UU,K'_p)|(\UU|\Oa')=0$ for all $q\ge1$ and $p\ge0$.
}
	
\biz{Proof.}{The proof is by simple and standard homological algebra based on
	a double complex.
	 Recall the complex $(\eBA)$ and note that its differential $d$ can be
	extended to (alternating) cochains componentwise, and that this extension, also called $d$,
	commutes with the \v Cech differential $\da$.
	 Given a cocycle $f\in Z^q(\UU,K'_p)$, $q\ge1$, $p\ge0$, of the finite covering
	$\UU$, we need to find a cochain $g\in C^{q-1}(\UU|\Oa',K'_p)$ with $\da g=f|\Oa'$.
	 To do that we determine cochains $\fii_i\in C^{q+i-1}(\UU,\OO^{L_{p+i}})$
	and $\psi_i\in C^{q+i-2}(\UU|\Oa',\OO^{L_{p+i}})$ for $i\ge1$ \st\
	$f=d\fii_1$, $\da\fii_i=d\fii_{i+1}$, and $\fii_i|\Oa'=d\psi_{i+1}+\da\psi_i$
	for $i\ge1$.

	 As \p\tBC(b) shows that $(\eBA)$ is $d$-exact on the global level 
	over the bodies of the simplices of $\UU$
	we can one after another find $\fii_1,\fii_2,\fii_3,\ldots$ \st\
	$f=d\fii_1$, $\da\fii_1=d\fii_2$, $\da\fii_2=d\fii_3$, $\ldots$.
	 Since $\fii_{i+1}=0$ for $i$ large enough as
	$\fii_{i+1}\in C^{q+i}(\UU,\OO^{L_{p+i+1}})$ and this group is zero
	for $i$ large enough because $\UU$ is finite, we see that
	$\da\fii_i=d\fii_{i+1}=0$, i.e., $\fii_i|\Oa'\in C^{q+i-1}(\UU|\Oa', \OO^{L_{p+i}})$
	is a $\da$-cocycle, and hence $\fii_i|\Oa'=\da\psi_i$, where
	$\psi_i\in C^{q+i}(\UU|\Oa',\OO^{L_{p+i}})$.
	 Let $\psi_j=0$ for $j>i$, and determine $\psi_{i-1},\psi_{i-2},\ldots,\psi_1$
	one after another.
	 
	 As $\da\fii_{i-1}|\Oa'=d\fii_i|\Oa'=\da d\psi_i$ we find that $\da(\fii_{i-1}|\Oa'-d\psi_i)=0$,
	i.e., $\fii_{i-1}|\Oa'=d\psi_i+\da\psi_{i-1}$, etc., $\fii_1|\Oa'=d\psi_2+\da\psi_1$.
	 Then as $f|\Oa'=d\fii_1|\Oa'=\da d\psi_1$ letting $g=d\psi_1$ will do.
	 The proof of \p\tBE\ is complete.
}

\alcim{\sC. EXHAUSTION.}

	 This section describes a way to exhaust a \pscx\ open subset $\Oa$ of a \bspc\ $X$
	that is convenient for proving vanishing results for sheaf cohomology over $\Oa$.
	 We follow here [\rLC, \S\,2].

	 We say that a \fn\ $\aa$, call their set $\AA'$, is an
	\kiemel{admissible radius \fn} on $\Oa$ if $\aa:\Oa\to(0,1)$ is \cts\
	and $\aa(x)<\dist(x,X\setminus\Oa)$ for $x\in\Oa$.
	 We say that a \fn\ $\aa$, call their set $\AA$, is an
	\kiemel{admissible Hartogs radius \fn} on $\Oa$ if $\aa\in\AA'$
	and $-\log\aa$ is \psh\ on $\Oa$.
	 Call $\AA$ \kiemel{cofinal} in $\AA'$ if for each $\aa\in\AA'$
	there is a $\ba\in\AA$ with $\ba(x)<\aa(x)$ for $x\in\Oa$.

\tetel{\p\tCA.}{Plurisubharmonic domination holds in\/ $\Oa$ if and only if $\AA$
	is cofinal in $\AA'$.
}

\biz{Proof.}{Write $\aa=e^{-u}\in\AA'$ and $\ba=e^{-\psi}\in\AA$.
	As \pshdom\ holds on $\Oa$ for $u$ \cts\ if and only if for
	$u$ locally upper bounded, the proof of \p\tCA\ is complete.
}

	 It will be often useful to look at coverings by balls
	$B_X(x,\aa(x))$, $x\in\Oa$, $\aa\in\AA'$ and shrink their
	radii to obtain a finer covering by balls $B_X(x,\ba(x))$,
	$x\in\Oa$, $\ba\in\AA$.

	 Let $e'_n$, $n\ge1$, be a \sbs\ in the \bspc\ $(X',\|\cdot\|')$,
	and similarly $e''_n$, $n\ge1$, in $(X'',\|\cdot\|'')$.
	 One can change the norms $\|\cdot\|',\|\cdot\|''$ to equivalent
	norms so that
$$
	\Big\|\sum_{i=m}^nx'_ie'_i\Big\|'\le\Big\|\sum_{i=M}^Nx'_ie'_i\Big\|'
\tag\eCA
 $$
 	for $0\le M\le m\le n\le N\le\infty$, $x'_i\in\CC$, and
	similarly for $X'',\|\cdot\|'',\{e''_i\}$.
	 Let $X=X'\times X''$ with norm $\|x\|=\max\{\|x'\|',\|x''\|''\}$
	on $x=(x',x'')\in X$, and \sbs\ $e_{2n-1}=e'_n$, $e_{2n}=e''_n$,
	for $n\ge1$, if both $X',X''$ are \idml.
	 If $k'=\dim(X')<\infty$, then let $e_n=e'_n$ for $n\le k'$,
	and $e_n=e''_{n-k'}$ for $n>k'$.
	 We assume that $X$ is \idml, since otherwise \t\tAB\ reduces to a
	well-known classical theorem.
	 Then $X,\|\cdot\|,\{e_i\}$ also satisfy the analog of $(\eCA)$.
	 Introduce the projections $\pi_N:X\to X$,
	$\pi_N\sum_{i=1}^\infty x_ie_i=\sum_{i=1}^Nx_ie_i$, 
	$x_i\in\CC$, $\pi_0=0$, $\pi_\infty=1$, $\ro_N=1-\pi_N$, and define for
	$\aa\in\AA$ and $N\ge0$ integer the sets
$$\eqalign{
	D_N\<\aa\>&=\{\xi\in\Oa\cap\pi_N X:(N+1)\aa(\xi)>1\}\cr
	\Oa_N\<\aa\>&=\{x\in\pi_N^{-1}D_N\<\aa\>:\|\ro_N x\|<\aa(\pi_N x)\}\cr
}.
\tag\eCB
 $$

	 These $\Oa_N\<\aa\>$ are \pscx\ open in $\Oa$, and they will serve
	to exhaust $\Oa$ as $N=0,1,2,\ldots$ varies.

\tetel{\p\tCB.}{Let $\aa\in\AA$, and suppose that \pshdom\ holds in\/ $\Oa$.
\vskip0pt  
	{\rm(a)} There is an $\aa'\in\AA$, $\aa'<\aa$, with\/
	$\Oa_n\<\aa'\>\subset\Oa_N\<\aa\>$ for all $N\ge n$.
	 So any $x_0\in\Oa$ has a \nbd\ contained in all but
	finitely many $\Oa_N\<\aa\>$.
\vskip0pt  
	{\rm(b)} There are $\ba,\ga\in\AA$, $\ga<\ba<\aa$, so that
	for all $N$ and $x\in\Oa_N\<\ga\>$
$$
	B_X(x,\ga(x))\subset\Oa_N\<\ba\>\cap\pi_N^{-1}B_X(\pi_N x,\ba(x))
	\subset B_X(x,\aa(x)),
\tag\eCC
 $$
\vskip0pt  
	{\rm(c)} and, additionally, for all $N$ there is a finite set of
	points $\xi_i=\pi_N\xi_i\in D_N\<\ga\>$ \st\ for each $x\in\Oa_N\<\ga\>$
	there is a $\xi_i$ with
$$
	B_X(x,\ga(x))\subset\Oa_N\<\ba\>\cap\pi_N^{-1}B_X(\xi_i,\ba(\xi_i)).
\tag\eCD
 $$
\vskip0pt
	{\rm(d)} There is an $\aa'\in\AA$, $\aa'<\aa$, with $\Oa_N\<\aa'\>\subset
	\Oa_N\<\aa\>\cap\Oa_{N+1}\<\aa\>$ for all $N\ge0$.
}

\biz{Proof.}{See [\rLC, Prop.\,2.1] and [\rLB, Prop.\,4.3] for (a) and (b),
	and [\rLC, Prop.\,2.3] for (d).
	We modify slightly Lempert's definition of $\ba,\ga$
	in his proof of (b) in [\rLC, Prop.\,2.1] so as to work also  for (c) here.

	To complete part (b) choose the functions $\ba$ and $\ga$ first in $\AA'$
	then in $\AA$ applying \pshdom\ in $\Oa$ so that $\ba<\aa/8$, $\ga<\ba/8$,
	$\aa(x)<2\aa(y)$ for $x,y\in B_X(z,2\ba(z))$, and
	$\ba(x)<2\ba(y)$ for $x,y\in B_X(z,2\ga(x))$.
	 Then (b) is verified as in the proof of [\rLC, Prop.\,2.1] arguing with the
	triangle inequality only.

	(c) Let $\epsz_0=\frac12\min_{D_N\<\ga\>}\ga$.
	As $\ga$ is strictly positive and \cts\ on the compact set $\overline{D_N\<\ga\>}$
	we see that $\epsz_0>0$.
	 Choose a finite $\epsz_0$-net $\{\xi_i\}$ in the totally bounded set
	$D_N\<\ga\>$, i.e., for each $\xi\in D_N\<\ga\>$ there is a $\xi_i$ with
	$\|\xi-\xi_i\|<\epsz_0$.
	 We claim that this choice of points $\{\xi_i\}$ will do.
	 Indeed, let $x\in\Oa_N\<\ga\>$ be any point, and pick a $\xi_i$
	with $\|\pi_Nx-\xi_i\|<\epsz_0$.
	 We need to show that if $\|y-x\|<\ga(x)$, then
	$y\in\Oa_N\<\ba\>$, and $\|\pi_Ny-\xi_i\|<\ba(\xi_i)$.
	 The already proved part (b) implies that $y\in\Oa_N\<\ba\>$, and we have
$$
\eqalign{
	\|\pi_Ny-\xi_i\|&\le\|\pi_N(y-x)\|+\|\pi_Nx-\xi_i\|
	\le
	\|y-x\|+\epsz_0\cr
	&<
	\ga(x)+\ga(\xi_i)
	<
	{\ts\frac18}\ba(x)+{\ts\frac18}\ba(\xi_i)\
	<
	{\ts\frac14}\ba(\pi_Nx)+{\ts\frac18}\ba(\xi_i)\cr
	&<
	{\ts\frac12}\ba(\xi_i)+{\ts\frac18}\ba(\xi_i)
	<
	\ba(\xi_i),\cr
}
 $$
 	where we used the above properties of $\ba$ and $\ga$ including the
	doubling inequality of $\ba$.
	 The proof of \p\tCB\ is complete.
}

	The meaning of \p\tCB(bc) is that certain refinement maps exist
	between certain open coverings.

\alcim{\sD. VANISHING.}

	This section completes the proof of \t\tAB.
	Resume the notation and hypotheses of \S\,\sA--\sC.
	For $\aa\in\AA$, $N\ge0$, put
$$
	\BB(\aa)=\{B_X(x,\aa(x)):x\in\Oa\},\qquad\BB_N(\aa)=\{B_X(x,\aa(x)):x\in\Oa_N\<\aa\>\}.
 $$

 	We say that a complex $(\eBA)$ satisfies
	\kiemel{condition~(\eDA)}
	(i.e., is tractable by our current methods)
	if there is an $\aa_0\in\AA$ so that 
	for any $\aa\in\AA$, $\aa<\aa_0$, there are $\ba,\ga\in\AA$ \st\
	\p\tCB\ holds for $\aa,\ba,\ga$ and the coverings $\BB(\aa)$,
	$\BB(\ga)$, 
	$\{\pi_N^{-1}B_X(\pi_Nx,\ba(x))\cap B_X(x,\aa(x)):x\in\Oa_N\<\ga\>\}$
	are Leray coverings of their respective unions for all $K'_p$, $p\ge1$.

	The following more natural condition below on a complex $(\eBA)$
	implies condition $(\eDA)$ above.

	We say that a complex $(\eBA)$ satisfies
	\kiemel{condition~(\eDB)}
	if there is an open covering $\UU$ of $\Oa$ by \pscx\ open subsets
	$U$ of $\Oa$ \st\ if $V$ is any \pscx\ open subset of any member $U\in\UU$,
	then $(\eBA)$ is exact over $V$ on the level of global sections.
	
	The reason for using the artificial looking condition $(\eDA)$ is that
	it is easily verifed a priori for our Koszul complex $(\eAA)$ while
	$(\eDB)$ not --- our \t\tAB\ is in fact equivalent to saying that
	$(\eAA)$ satisfies condition $(\eDB)$.
	 It seems unknown whether acyclicity holds over $\Oa$ for all
	complexes $(\eBA)$ in general.

\tetel{\p\tDA.}{The Koszul complex~$(\eAA)$ satisfies condition~$(\eDA)$.
}

\biz{Proof.}{We saw that there is an $\aa_0\in\AA'$ \st\ if $\aa\in\AA$,
	$\aa<\aa_0$, then $\BB(\aa)$ is a Leray covering, and we know from
	\p\tBA\ that $(\eAA)$ is acyclic on the level of global sections
	over any set 
	(e.g., $\pi_N^{-1}B_X(\pi_Nx,\ba(x))\cap B_X(x,\aa(x))$)
	that is the product of a bounded convex open set in $X'$ by
	a bounded convex open set in $X''$.
	 The proof of \p\tDA\ is complete.
	
}
	
\tetel{\p\tDB.}{Let $X$ be a \bspc\ with a \sbs, $\Oa\subset X$ \pscx\ open,
	$(\eBA)$ a complex, and suppose that \pshdom\ holds in any \pscx\ open
	subset of\/ $\Oa$, and that\/ $(\eBA)$ satisfies condition $(\eDA)$.
	Then for any $\aa\in\AA$ there is a $\ga\in\AA$ with $\ga<\aa$ and
	$H^q(\BB_N(\aa),K'_p)|\BB_N(\ga)=0$ for all $N\ge0$, $q\ge1$, and $p\ge0$.
}

\biz{Proof.}{We consider some open coverings and refinement maps of them.
	Let $\aa,\ba,\ga,\{\xi_i\}$ be as in \p\tCB(bc).
	Consider the open coverings
	$\BB_N(\aa)$,
	$\UU_N=\{U(x)=\pi_N^{-1}B_X(\pi_Nx,\ba(x))\cap B_X(x,\aa(x)):x\in\Oa_N\<\ga\>\}$,
	$\UU'_N=\{U(\xi_i):i\}$,
	$\BB_N(\ga)$, and their refinement maps
	$\UU_N\to\BB_N(\aa)$ given by
	$U(x)\mapsto B_X(x,\aa(x))$,
	$\UU'_N\to\UU_N$ given by
	$U(\xi_i)\mapsto U(\xi_i)$, and
	$\BB_N(\ga)\to\UU'_N|\Oa_N\<\ba\>$ given by
	$B_X(x,\ga(x))\mapsto\Oa_N\<\ba\>\cap U(\xi_i)$,
	where $\xi_i$ is the point assigned to $x\in\Oa_N\<\ga\>$ in
	\p\tCB(c).
	 Due to the inequalities $(\eCC)$ and $(\eCD)$ the above
	are indeed refinement maps and hence induce maps
$$\eqalign{
	H^q(\BB_N(\aa),&K'_p)\to
	H^q(\UU_N,K'_p)\to
	H^q(\UU'_N,K'_p)\to\cr
	&H^q(\UU'_N,K'_p)|(\UU'_N|\Oa_N\<\ba\>)\to
	H^q(\BB_N(\ga),K'_p)\cr
}
\tag\eDC
 $$
 	in cohomology for $q\ge1$, and $p\ge0$.
	 We see via condition $(\eDA)$ that \p\tBE\ applies, and thus
	the fourth group in $(\eDC)$ is zero, and together with it so is
	the composite map $H^q(\BB_N(\aa),K'_p)\to H^q(\BB_N(\ga),K'_p)$ in $(\eDC)$.
	 The proof of \p\tDB\ is complete.
}

	 The upshot of \p\tDB\ is that the refinement map $\BB_N(\ga)\to\BB_N(\aa)$
	factors through a finite covering.

\tetel{\t\tDC.}{Let $X,\Oa,(\eBA)$ be as in \p\tDB.
\vskip0pt
	{\rm(a)} $H^q(\Oa,K'_p)=0$ for all $q\ge1$ and $p\ge0$.
\vskip0pt
	{\rm(b)} The sequence\/ $(\eBA)$ is exact over\/ $\Oa$
	on the level of global sections.
}

\biz{Proof.}{(a) Suppose first that $q\ge2$.
	Let $f\in H^q(\Oa,K'_p)$ be a cohomology class.
	 We saw earlier that due to \pshdom\ in $\Oa$ there is an $\aa$ with
	$10\aa\in\AA$ so that $f$ can be represented as a cocycle
	$f\in Z^q(\BB(\aa),K'_p)$.
	 Choose $\ga$ as in \p\tDB.
	 \p\tDB\ yields a $g_N\in C^{q-1}(\BB_N(\ga),K'_p)$ with $\da g_N=f|\BB_N(\ga)$.
	 We can extend the cochain $g_N$ to a cochain $g_N\in C^{q-1}(\BB(\ga),K'_p)$
	simply by defining $g_N$ to be zero over simplices
	$\bigcap_{i=1}^q B_X(x_i,\ga(x_i))$ if at least one vertex $x_i\not\in\Oa_N\<\ga\>$.
	 \p\tDB\ gives a $\ga'\in\AA$, $\ga'<\ga$, with $H^{q-1}(\BB_{N-1}(\ga),K'_p)|\BB_{N-1}(\ga')=0$
	for all $N\ge1$, and \p\tCB(d) a $\ga''\in\AA$, $\ga''<\ga'$ with
	$\Oa_{N-1}\<\ga''\>\subset\Oa_{N-1}\<\ga'\>\cap\Oa_N\<\ga'\>$ for all $N\ge1$.
	 So similarly by extending a $(q-2)$-cochain there is an 
	$h_N\in C^{q-2}(\BB(\ga''),K'_p)$ with
	$(g_N-g_{N-1})|\BB_{N-1}(\ga'')=\da h_N|\BB_{N-1}(\ga'')$.
	 Letting $g'_N=g_N|\BB(\ga'')-\sum_{n=1}^N\da h_n
	\in C^{q-1}(\BB(\ga''),K'_p)$
	\p\tCB(a) implies as $g'_N|\BB_{N-1}(\ga'')=g'_{N-1}|\BB_{N-1}(\ga'')$ that
	$g'_N$ converges as $N\to\infty$ in a quasistationary manner to a
	$g\in C^{q-1}(\BB(\ga''),K'_p)$ with $\da g=f|\BB(\ga'')$.
	 Thus $f$ equals zero in $H^q(\Oa,K'_p)$ for $q\ge2$ and $p\ge0$.
         Part (a) follows then from \p\tBC(c) and (b) from (a) via \p\tBC(b).
	 The proof of \t\tDC\ is complete.

}

\biz{Proof of \t\tAB.}{Propositions \tBB\ and \tDA\ exhibit \t\tAB\ as a special case of \t\tDC.
	The proof of \t\tAB\ is complete.
}

\alcim{\sE. AN ANALOG OF HEFER'S LEMMA.} 

	 In this section we extend \t\tAB\ and draw some corollaries
	from it.
	 
	 While it seems far from being currently proved, it is reasonable to hope
	that \pshdom\ holds in every \pscx\ open subset $\Oa$ of any \bspc\ $X$ that
	is a direct summand of a \bspc\ $Y$, which has a \sbs\ (i.e., if $X$ has the
	bounded approximation property, fondly called BAP). It is certainly a good
	question to ask.
	 In this spirit we can relax the hypotheses of \t\tAB\ as follows.

\tetel{\t\tEA.}{Let $X',X'',Z$ be \bspc{}s, $\Oa\subset X=X'\times X''$ \pscx\ open,
	$I$ the sheaf of germs of \holo\ \fns\ $\Oa\to Z$ that vanish on $X'$.
	 Suppose that $X$ has a \sbs, and that \pshdom\ holds in every \pscx\ open
	subset of\/ $\tilde\Oa=\Oa\times X'\times X''$. Then
{\vskip0pt\rm (a)}
	the Koszul complex\/ $(\eAA)$ is exact on the germ level and on the level
	of global sections over\/ $\Oa$, and
{\vskip0pt\rm (b)}
	the $K_p$ are acyclic over\/ $\Oa$: $H^q(\Oa,K_p)=0$ for all $q\ge1$ and $p\ge0$.
}

\biz{Proof.}{We consider a Koszul complex $(\eEA)$ as in $(\eAA)$ over a bigger \bspc.
	Let $\tilde x=(x',x'',y',y'')\in\tilde X=X'\times X''\times X'\times X''=
	\tilde X'\times\tilde X''\ni(\tilde x'=(x',y''),\tilde x''=(x'',y'))$,
	$\tilde\La_p$ the \bspc\ of all \cts\ complex $p$-linear alternating maps
	$\tilde X''\to Z$ for $p\ge0$, $\tilde\La_0=\tilde\La_{-1}=Z$, $\OO^{\tilde\La_p}\to\tilde X$
	the sheaf of germs of \holo\ \fns\ $\tilde X\to\tilde\La_p$, $\tilde E$ the Euler
	vector field on $\tilde X''$ defined by $\tilde E(\tilde x'')=\tilde x''=(x'',y')$,
	$i_{\tilde E}$ the inner derivation determined by $\tilde E$, i.e.,
	$i_{\tilde E}$ is the contraction of $p$-forms with $\tilde E$:
	if $f$ is a local section of $\OO^{\tilde\La_p}$, then let $i_{\tilde E}f$ be
	the local section of $\OO^{\tilde\La_{p-1}}$ given for $p\ge1$ by
	$(i_{\tilde E}f)(x',x'',y',y'')((\xi''_1,\eta'_1),\ldots,(\xi''_{p-1},\eta'_{p-1}))=
	f(x',x'',y',y'')((x'',y'),(\xi''_1,\eta'_1),\ldots,(\xi''_{p-1},\eta'_{p-1}))$,
	and for $p\ge0$ by the formula $(i_{\tilde E}f)(x',x'',y',y'')=f(x',0,0,y'')$.
	 Let $\tilde I$ be the subsheaf of $\OO^Z\to\tilde X$ of sections that vanish on
	$\tilde X'$.
	 We consider the Koszul complex
$$
	\ldots\to\OO^{\tilde\La_p}\to\OO^{\tilde\La_{p-1}}\to\ldots\to\OO^{\tilde\La_1}\to\tilde I\to0
\tag\eEA
 $$
 	of analytic sheaves over $\tilde X$, where each map is $i_{\tilde E}$.
	 Let $\tilde K_p$, $p\ge0$, be the corresponding sequence of kernel sheaves:
	$\tilde K_p(U)=\{f\in\OO(U,\tilde\La_p):i_{\tilde E}f=0{\text{\/ on\ }}U\}$,
	$U\subset\tilde X$ open; $\tilde K_0=\tilde I$.
	 Let $\tilde\Oa=\{(x',x'',y',y'')\in\tilde X:(x',x'')\in\Oa\}=\Oa\times X'\times X''$.

	 \t\tAB\ applies to $(\eEA)$ over $\tilde\Oa$ since $\tilde X'\cong\tilde X''=X=X'\times X''$
	have a \sbs\ by assumption.
	 Now \t\tEA\ follows easily by considering the extension of forms $f\in\OO(U,\La_p)$
	to forms $\tilde f\in\OO(\tilde U,\tilde\La_p)$, where $U\subset\Oa$ open,
	$\tilde U=U\times X'\times X''$, $p\ge0$, and the restriction $f$ of forms
	$\tilde f$ to $\Oa$, defined for $f\mapsto\tilde f$ by
	$\tilde f(x',x'',y',y'')((\xi''_1,\eta'_1),\ldots,(\xi''_p,\eta'_p))=
	f(x',x'')(\xi''_1,\ldots,\xi''_p)$, and for $\tilde f\mapsto f$ by
	$f(x',x'')(\xi''_1,\ldots,\xi''_p)=\tilde f(x',x'',0,0)((\xi''_1,0),\ldots,(\xi''_p,0))$.
	 As these simple maps intertwine $i_E$ and $i_{\tilde E}$, the proof of \t\tEA\ is complete.
}

	 Similarly one could replace in \t\tEA\ the assumption that `$X$ have a \sbs' by `$X$ be a direct
	 summand of a \bspc\ $Y$ with a \sbs.'

	 The next item is an \idml\ analog of the classical Hefer lemma (without bounds).

\tetel{\t\tEB.}{Let $X,Z$ be \bspc{}s, $\Oa\subset X$ \pscx\ open, and $f\in\OO(\Oa\times\Oa,Z)$.
	If $X$ has a \sbs\ and \pshdom\ holds in every \pscx\ open subset of\/ $\Oa\times\Oa$,
	and $f(x,x)=0$ for $x\in\Oa$, then there is a $g\in\OO(\Oa\times\Oa,\Hom(X,Z))$ \st\
	$f(x,y)=g(x,y)(x-y)$ for $x,y\in\Oa$.
}

\biz{Proof.}{Writing $x'=\frac12(x+y)$, $x''=\frac12(x-y)$, $x'+x''=x$, $x'-x''=y$
	a direct decomposition of $(x,y)\in X\times X\cong X'\times X''\ni(x',x'')$ is obtained
	with $X'\cong X''\cong X$.
	 As $f=0$ for $x''=0$ an application of \t\tAB(a) completes the proof of \t\tEB.
}

	 The above form of Hefer's lemma can be applied to give an algebraic definition
	of the Fr\'echet differential $df\in\OO(\Oa,X^*)$ of a \fn\ $f\in\OO(\Oa)$
	defined on a \pscx\ open subset $\Oa$ of a \bspc\ $X$, which is just like the
	classical case of polynomials in finitely many variables.

\tetel{\t\tEC.}{Let $X,Z$ be \bspc{}s, $\Oa\subset X$ \pscx\ open, and $f\in\OO(\Oa,Z)$.
	If $X$ has a \sbs\ and \pshdom\ holds in every \pscx\ open subset of\/ $\Oa\times\Oa$,
	then the \fn\ $f(x)-f(y)$ can be written as $f(x)-f(y)=g(x,y)(x-y)$ for $x,y\in\Oa$,
	where $g\in\OO(\Oa\times\Oa,\Hom(X,Z))$, and
	$df(x)\xi=g(x,x)\xi$ for $x\in\Oa$, $\xi\in X$.
}

\biz{Proof.}{\t\tEB\ provides such a $g$.
	Taking an $x$ partial derivative of the identity in \t\tEB\
	in the $\xi$ direction and letting $y=x$ complete the proof of \t\tEC.
}

	 \t\tEB\ can also be formulated for arbitrary \bvbdl{}s instead of just
	trivial ones.

\tetel{\t\tED.}{Let $X$ be a \bspc, $\Oa\subset X$ \pscx\ open, $E\to\Oa\times\Oa$
	a \holo\ \bvbdl, and $f\in\OO(\Oa\times\Oa,E)$ a \holo\ section.
	 If $X$ has a \sbs\ and \pshdom\ holds in every \pscx\ open subset of\/
	$\Oa\times\Oa$, and $f(x,x)=0$ for all $x\in\Oa$, then there is a
	$g\in\OO(\Oa\times\Oa,\Hom(X,E))$ \st\ $f(x,y)=g(x,y)(x-y)$ for $x,y\in\Oa$.
}

\biz{Proof.}{\t1.3(b) in [\rPA] gives a \bspc\ $Z_1$, an embedding $I\in\OO(\Oa\times\Oa,\Hom(E,Z_1))$,
	and a projection $P\in\OO(\Oa\times\Oa,\End(Z_1))$ that $I(x,y)(E_{x,y})=P(x,y)(Z_1)$
	for all $x,y\in\Oa$.
	 Look at the \fn\ $f'\in\OO(\Oa\times\Oa,Z_1)$ defined by $f'(x,y)=I(x,y)f(x,y)$,
	which vanishes for $x=y$ in $\Oa$.
	 \t\tEB\ gives a $g'\in\OO(\Oa\times\Oa,\Hom(X,Z_1))$ with $f'(x,y)=g'(x,y)(x-y)$.
	 Define $g\in\OO(\Oa\times\Oa,\Hom(X,E))$ by
	$g(x,y)=I(x,y)^{-1}P(x,y)g'(x,y)$.
	 As $g(x,y)(x-y)=I(x,y)^{-1}P(x,y)I(x,y)f(x,y)=I(x,y)^{-1}I(x,y)f(x,y)=f(x,y)$,
	since $P(x,y)I(x,y)$ is the identity, the proof of \t\tED\ is complete.
}

	 \t\tEA, just like \t\tEB, has a version with \bvbdl{}s.

	 Let $X',X''$ be \bspc{}s, $X=X'\times X''$, $\Oa\subset X$ \pscx\ open,
	$F\to\Oa$ a \holo\ \bvbdl, $\La_p(F)\to\Oa$ the \bvbdl\ of all the \bspc{}s
	$\La_p(F_x)$, $x\in\Oa$, of all \cts\ complex $p$-linear alternating maps
	from $X''$ to the fiber $F_x$ of $F$ over the point $x$, $p\ge0$;
	$\La_0(F_x)=F_x$, and $\OO^{\La_p(F)}\to\Oa$ the sheaf of germs of \holo\
	sections $\La_p(F)\to\Oa$.
	 Let $E$ be the Euler vector field on $X''$ defined by $E(x'')=x''$, $i_E$
	the inner derivation determined by $E$, i.e., $i_E$ is the contraction of
	$p$-forms in $\La_p(F)$ with $E$: if $f$ is a local section of $\OO^{\La_p(F)}$,
	then let $i_Ef$ be the local section of $\OO^{\La_{p-1}(F)}$ given for $p\ge1$
	by $(i_Ef)(x',x'')(\xi''_1,\ldots,\xi''_{p-1})=f(x',x'')(x'',\xi''_1,\ldots,\xi''_{p-1})$.
	 Let $I$ be the subsheaf of $\OO^F\to\Oa$ of all sections that vanish on 
	$\Oa\cap (X'\times\{0\})$.
	 We consider the Koszul complex
$$
	\ldots\to\OO^{\La_p(F)}\to\OO^{\La_{p-1}(F)}\to\ldots\to\OO^{\La_1(F)}\to I\to0
\tag\eEB
 $$
 	of analytic sheaves over $\Oa$, where each map except the rightmost is $i_E$.
	 Let $K_p$, $p\ge0$, be the corresponding sequence of kernel sheaves:
	$K_p(U)=\{f\in\OO(U,\La_p(F)):i_Ef=0{\text{\/ on\ }}U\}$, $U\subset\Oa$ open, $p\ge1$;
	$K_0=I$.

\tetel{\t\tEE.}{With the above notation, suppose that $X$ has a \sbs, and that \pshdom\
	holds in every \pscx\ open subset of\/ $\Oa\times X$.
	 Then
{\vskip0pt\rm (a)}
	the Koszul complex\/ $(\eEB)$ is exact on the germ level and on the level
	of global sections over\/ $\Oa$, and
{\vskip0pt\rm (b)}
	the $K_p$ are acyclic over\/ $\Oa$: $H^q(\Oa,K_p)=0$ for all $q\ge1$ and $p\ge0$.
}

\biz{Proof.}{Just like in the proof of \t\tED, rely on \t1.3(b) in [\rPA] to exhibit
	the \bvbdl\ $F$ as a direct summand of a trivial \bvbdl\ $\Oa\times Z_1$, and thus
	the complex $(\eEB)$ as a direct summand of a complex $(\eAA)$.
	 An application of \t\tEA\ completes the proof of \t\tEE.
}

	 We conclude by remarking that the Koszul complex $(\eAA)$ of the ideal
	sheaf of the origin in a \bspc\ $X$ is also useful in connection with
	the projectivization $P(X)$ of $X$, and that Hefer's lemma as above
	helps us understand the universal derivation of the algebra $\OO(\Oa)$
	for $\Oa\subset X$ \pscx\ open; see [\rNW] for the latter if $X$ is \fdml.

	{\vekony Acknowledgements.} The work on this paper started at the Riverside
	campus of the University of California, continued at its San Diego
	campus, and was completed at Georgia State University.
	 The author is grateful to these institutions and to 
	Professors Bun Wong of UCR, Peter Ebenfelt of
	UCSD, and Mih\'aly Bakonyi of GSU, whose support has been essential to him.

\vskip0.30truein
\centerline{\scVIII References}
\vskip0.20truein
\baselineskip=11pt
\parskip=7pt
\frenchspacing
{\rmVIII

	[\rLA] Lempert,~L.,
	{\itVIII
	Plurisubharmonic domination},
	J. Amer. Math. Soc.,
	{\bfVIII 17}
	(2004),
	361--372.

	[\rLB] \vonal,
	{\itVIII
	Vanishing cohomology for holomorphic vector bundles in a Banach setting},
	Asian J. Math.,
	{\bfVIII 8}
	(2004),
	65--85.

    	[\rLC] \vonal,
    	{\itVIII
	Acyclic sheaves in Banach spaces},
	Contemporary Math.,
	{\bfVIII 368}
	(2005),
	313--320.

	[\rNW] Neeb,~K.-H., Wagemann,~F.,
	{\itVIII
	The universal central extension of the holomorphic current algebra},
	Manuscripta Math.,
	{\bfVIII 112}
	(2003),
	441--458.

	[\rPA] Patyi,~I., 
	{\itVIII
	On holomorphic Banach vector bundles over Banach spaces},
	manuscript, (2005).

	[\rPB] \vonal, 
	{\itVIII
	Analytic cohomology in a Banach space},
	in preparation.

\vskip0.20truein
\centerline{\vastag*~***~*}
\vskip0.15truein
{\scVIII
	Imre Patyi,
	Department of Mathematics and Statistics,
	Georgia State University,
	Atlanta, GA 30303-3083, USA,
	{\ttVIII ipatyi\@gsu.edu}}
\bye